\documentclass[11pt, reqno]{article}
\usepackage{amsfonts}
\usepackage{amsmath}
\usepackage{amsthm}
\usepackage{amssymb}
\def\b{\mathbb }

\def\phi{\varphi }

\def\reel{\b R}

\def\comp{\b C}
\def\ganz{\b Z}
\def\hat{\widehat}

\swapnumbers
\theoremstyle{plain}
\newtheorem{theorem}{Theorem}[section]
\newtheorem{corollary}[theorem]{Corollary}
\newtheorem{lemma}[theorem]{Lemma}
\newtheorem{proposition}[theorem]{Proposition}
\theoremstyle{definition}

\newtheorem{remarks}[theorem]{Remarks}
\newtheorem{examples}[theorem]{Examples}

\newtheorem{Duop}[theorem]{Dunkl operators}
\newtheorem{refl}[theorem]{Basic notions}
\newtheorem{Dutr}[theorem]{The Dunkl transform}

\newtheorem{Lap}[theorem]{The generalized Laplacian}
\newtheorem{heat}[theorem]{Generalized heat kernels}
\newtheorem{Gau}[theorem]{$k$-Gaussian semigroups}

\newtheorem{Mom}[theorem]{Modified moment functions}
\newtheorem{Mod}[theorem]{Modified moments of $k$-Gaussian measures}

\newtheorem{app}[theorem]{ Appell characters}

\newtheorem{appco}[theorem]{ Appell cocharacters}

\numberwithin{equation}{section}

\begin{document}
\title{Biorthogonal polynomials associated with  reflection groups\\
and a formula of Macdonald}
\author{
Margit R\"osler\thanks{Parts of this paper were written  while the  author 
held a Forschungsstipendium of the DFG at the University of Virginia, 
Charlottesville, USA.} \\
 Zentrum Mathematik, Technische Universit\"at M\"unchen\\
 Arcisstr. 21, D-80290 M\"unchen, Germany\\
e-mail: roesler@mathematik.tu-muenchen.de\\
and
 Michael Voit\thanks{Parts of this paper were written  while the  author
 was a visiting lecturer
at the University of Virginia, 
Charlottesville, USA. }\\
 Mathematisches Institut,  Universit\"at T\"ubingen\\
Auf der Morgenstelle 10, 72076  T\"ubingen, Germany\\
e-mail: voit@uni-tuebingen.de\\
}

\date{}

\maketitle

\begin{abstract}  Dunkl operators are differential-difference operators on $\reel^N$
which  generalize partial derivatives.
  They lead to generalizations of Laplace operators,
  Fourier transforms,
 heat semigroups, Hermite polynomials, and so on.
 In this paper we  introduce 
two systems of biorthogonal polynomials with
 respect to  Dunkl's Gaussian distributions in
 a quite canonical way.  These systems, called   Appell systems,
 admit many properties known from classical Hermite
 polynomials, and turn out to be useful for the analysis of Dunkl's
 Gaussian distributions. In particular, these polynomials lead to a
 new proof of a generalized formula of Macdonald  due to
 Dunkl. The ideas for this paper are taken from recent works on non-Gaussian
 white noise analysis and from  the umbral calculus. 
\end{abstract}

\smallskip\smallskip
\section{Introduction}

 Dunkl operators are differential-difference operators on $\reel^N$
 related to finite reflection groups.
They can be regarded as a generalization of partial derivatives and play
 a major role in the description of
Calogero-Moser-Sutherland models of quantum many-body systems on the line. 
  Dunkl operators lead to generalizations of exponential functions,
 Fourier transforms, Laplace operators, and Gaussian distributions on
 $\reel^N$.
The corresponding basic theory is developed in [D1-3, dJ, O, R2-3] and
 will 
be briefly reviewed in Section 2 below (with references
 to proofs). A more detailed approach to the Dunkl theory will be
 given in [R-V]. 

In this paper, we study  systems of biorthogonal polynomials
with respect to generalized Gaussian distributions of the form
$w(x)\,\text{exp}\,\{-c|x|^2\}$, where $w$ is a homogeneous weight function with
certain reflection group symmetries. Generalized Gaussians of such
type were first considered by Mcdonald and Mehta (see \cite{Me}) and
replace the classical Gaussian in the Dunkl theory.  
We mention that systems of orthogonal polynomials related to Dunkl's  Gaussian
distributions on $\reel^N$,  
called generalized Hermite polynomials, have been studied in various
contexts during the last years, see e.g [B-F, vD,
R2, Ros] and references therein. In the case of the
symmetric
group,  the  generalized Hermite polynomials play a role in quantum
many body systems of 
Calogero-Moser-Sutherland type and are closely related to  
Jack polynomials. 

The  biorthogonal systems of this paper 
are introduced in a quite canonical way; 
the method is as follows: In Section 3 we define so-called
modified moment functions $m_\nu$ ($\nu\in\ganz_+^N$) which form
generalizations of the monomials $x^\nu:= x_1^{\nu_1}\cdots
x_N^{\nu_N}$ on $\reel^N$ such that each $m_\nu$ is a homogeneous
polynomial of degree $|\nu|=\nu_1+\ldots+\nu_N$. These functions lead
to generalized moments of Dunkl's Gaussian distributions and are very
closely related to
the classical moments of the classical Gaussian distributions on
$\reel^N$. Based on these modified moments, we introduce two systems
$R_\nu(t,x)$ and $S_\nu(t,x)$ ($\nu\in\ganz_+^N$, $x\in\reel^N$,
$t\in\reel$)
of polynomials in $N+1$ variables via two generating functions involving
Dunkl's kernel $K$ and Dunkl's Gaussian distribution. We show that
for all $t>0$, the systems $(R_\nu(t,.))$ and 
 $(S_\nu(t,.))$ form biorthogonal polynomials with 
 respect to
Dunkl's Gaussian distribution with variance parameter $t$. We also
show that $R_\nu(t,x)=e^{-t\Delta_k}m_\nu(x)$ and
$S_\nu(t,x)=e^{-t\Delta_k} x^\nu$ hold where $\Delta_k$ denotes Dunkl's
Laplacian. At the end of this paper, we present some applications of
these results including a  Rodriguez-type formula and a new proof of a
generalized Macdonald-formula given originally in [D2]. We point out
that the methods of this paper can be extended to other distributions
and that there exist natural applications in martingale theory 
(which are well-known for classical Hermite polynomials and Brownian
motions); for details see [R-V]. The main purpose of this paper is
to give a short introduction to the Appell systems
 $(R_\nu(t,.))_{\nu\in\ganz_+^N}$ and 
 $(S_\nu(t,.))_{\nu\in\ganz_+^N}$ without a deeper probabilistic
 background. 

\section{Dunkl operators and Dunkl transform}

In this section we collect some basic   notations and  
  facts from the Dunkl theory which will be 
 important later on. General references here are  [D1, D2, D3, dJ, O];
 for a background on  reflection groups and root systems see, for instance, [Hu].

\begin{refl} For $\alpha \in \b R^N\setminus\{0\}$,  let $\sigma_\alpha$ be  the reflection 
 on the hyperplane $H_\alpha\subset\reel^N$ 
orthogonal to $\alpha$, i.e.,
$ \sigma_\alpha (x) = x- (2\langle\alpha,x\rangle/{|\alpha|^2})\,
\alpha,$
where  $\langle.,.\rangle$ is the Euclidean scalar product on $\b R^N$ and 
$|x|:= \sqrt{\langle x,x\rangle}$.

A finite set $R\subset\reel^N\setminus\{0\}$ is called a root system
 if $R\,\cap\, \reel\cdot\alpha=\{\pm \alpha\}$ and $\sigma_\alpha R=R\,$
 for all $\alpha\in R$. For a given
root system $R$ the reflections $\sigma_\alpha$ ($\alpha\in R$)
 generate a group $W,$  
the reflection group associated with $R$. This group  is finite,
 and all reflections in $W$ correspond to suitable pairs of roots; see [Hu].
Now fix $\beta\in\reel^N\setminus\bigcup_{\alpha\in R} H_\alpha$ and
  a  positive subsystem
 $R_+=\{\alpha\in R:\>\langle \alpha,\beta\rangle>0\}$; then 
 for each $\alpha\in R$ either $\alpha\in R_+$ or  $-\alpha\in R_+$.
We assume from now on with no loss of generality that the root
 system $R$ is normalized, i.e,  that $|\alpha|=\sqrt 2$ for all $\alpha\in R$.

A multiplicity function $k$ on  a root system $R$ is defined as a
 $W$-invariant function $k:R\to{\b C}$.
If one regards $k$ as function on the corresponding reflections, 
this $W$-invariance just means that $k$ is
 constant on the conjugacy classes of  reflections in $W.$ For abbreviation, we  introduce
\begin{equation}
\textstyle \gamma:=\gamma(k):=\sum_{\alpha\in R_+}k(\alpha)
\quad\quad{\rm and}\quad\quad
c_k:=\Bigl( \int_{\reel^N} e^{-|x|^2}\> w_k(x)\> dx\Bigr)^{-1}
\end{equation}
where $w_k$ is the  weight function 
\begin{equation}
 w_k(x) = \prod_{\alpha\in R_+}
 |\langle\alpha,x\rangle|^{2k(\alpha)}\,.
\end{equation}
Notice that $w_k$ is $W$-invariant  and homogeneous of degree $2\gamma$.
We shall use the following further abbreviations:  
 $\,\mathcal P=\b C\,[\b R^N]$ denotes the algebra of polynomial 
functions on $\b R^N$ and $\mathcal P_n \>\, (n\in \b Z_+)$ 
the subspace of  homogeneous polynomials of degree $n$. We use the
 standard multi-index notations, i.e. for multi-indices $\nu,\,\rho\in
 \b Z_+^N$ we write
$$|\nu|:= \nu_1+\ldots + \nu_N\,,\quad\quad \nu!:=\nu_1!\cdot \nu_2!\cdots \nu_N!\,,\quad\quad
\binom{\nu}{\rho} := \binom{\nu_1}{\rho_1} \binom{\nu_2}{\rho_2} \ldots
 \binom{\nu_N}{\rho_N}, $$
as well as
$$x^\nu:=x_1^{\nu_1}\cdots x_N^{\nu_N} \quad\text{ and }\quad 
A^\nu:=A_1^{\nu_1}\cdots A_N^{\nu_N}$$
for $x\in \b R^N$ and any family $A=(A_1,\ldots,A_N)$ of commuting
operators on $\mathcal P$. Finally,  the  partial ordering $\leq$ on $\b
Z_+^N$ is defined by $\,\rho\leq\nu :\Longleftrightarrow \rho_i \leq \nu_i$
for $i=1,\ldots, N.$ 
\end{refl}

\begin{Duop} The Dunkl operators $T_i \>\>(i=1,\ldots,N)$ on $\b R^N$ 
associated 
with the finite reflection group $W$ and multiplicity function $k$ are  given by 
\begin{equation}
 T_if(x) := \partial_if(x) + \sum_{\alpha\in R_+} k(\alpha)\,\alpha_i\cdot 
\frac{f(x)-f(\sigma_\alpha x)}{\langle\alpha,x\rangle}\,,\>\quad f\in 
C^1(\b R^N);
\end{equation}
here $\partial_i$ denotes the $i$-th partial derivative. In case $k=0$, the 
$T_i$ reduce to the corresponding partial derivatives. In this paper, we always
 assume  that $k\geq 0$ (i.e. all values of $k$ are non-negative), 
though several results  may be extended to larger ranges of $k$.  
The most important basic properties of the $T_i$  are 
as follows (see \cite{Du1}):
\begin{enumerate}
\item[(1)] The set $\{T_i\}$ generates a commutative algebra of 
differential-difference operators on $\mathcal P$.
\item[(2)] Each $T_i$ is homogeneous of degree $-1$ on $\mathcal P$, i.e.,
 $T_i \,p \in \mathcal P_{n-1}$ for $p\in\mathcal P_n.$ 
\item[(3)] (Product rule:) $ T_i(fg)= (T_i f)g + f(T_i g)$ for
$i=1,\ldots,N$ and  $f,g\in C^1(\b R^N)$ such that $g$ is 
 $W$-invariant.
\end{enumerate}
\end{Duop}

A major tool in this paper is a generalized exponential kernel $K(x,y)$ on 
$\b R^N \times \b R^N$, which replaces the usual exponential function  
$e^{\langle x,y\rangle}$. It  was  introduced in \cite{Du2} by means of 
an  intertwining isomorphism $V$ of $\cal P$ which is characterized by
the properties
\begin{equation}
V({\cal P}_n)={\cal P}_n,\quad V|_{{\cal P}_0}=id,
\quad{\rm and}\quad T_iV=V\partial_i \quad(i=1,\ldots,N).
\end{equation}
Let $B=\{x\in\b R^N: |x|\leq 1\}.$ Then $V$ extends to a bounded
linear operator on the algebra
\[ A:= \{f: B\to\b C: f=\sum_{n=0}^\infty f_n \>\>\>\text{ with }\,\, f_n\in
\mathcal P_n \>\>\>\text{and }\> \|f\|_A := \sum_{n=0}^\infty
\|f_n\|_{\infty,B} \,<\infty\}\]
by $\,V\bigl(\sum_{n=0}^\infty f_n\bigr) := \sum_{n=0}^\infty Vf_n\,.$ The Dunkl
kernel $K$ is now defined by  
\begin{equation}\label{(K)}
K(x,y):= V(e^{\langle.,\,y\rangle})(x) \quad\, (x,y\in \b R^N).
\end{equation}
 $K$  has a holomorphic extension 
to $\b C^N\times \b C^N$ and is symmetric in its arguments. We also
note that for $y\in\reel^N$, the function $ x\mapsto K(x,y)$  may be 
 characterized  as the unique analytic solution of 
$\, T_i f = y_i f\,\, (i=1,\ldots,N)$  
with  $f(0)=1$; see [O].

\begin{examples}
\begin{enumerate}
\item[\rm{(1)}]  If $k=0$, then $K(z,w)=e^{\langle z,w\rangle }$ 
for all $z,w\in {\b C}^N.$
\item[\rm{(2)}] If $N=1$ and   $W=\b Z_2$, then the  multiplicity function 
is a single parameter $k\geq 0$, and the intertwining operator is given
explicitly by
\[ V_kf(x) = \frac{\Gamma(k+1/2)}{\Gamma(1/2)\,\Gamma(k)}\int_{-1}^1
f(xt) (1-t)^{k-1}(1+t)^k\,dt;\]
see \cite{Du2}.  The associated Dunkl  kernel can be written as
\[ K(z,w) = j_{k-1/2}(izw) + \frac{zw}{2k +1}\,j_{k+1/2}(izw),
\quad z,w\in \b C,\]
where for $\alpha\geq -1/2$, $j_\alpha$ is the normalized spherical 
Bessel function 
\[ j_\alpha(z) = 2^\alpha\Gamma(\alpha +1) \frac{J_\alpha(z)}{z^\alpha}\,=\, 
\Gamma(\alpha +1) \cdot\sum_{n=0}^\infty \frac{(-1)^n(z/2)^{2n}}{n! \,
\Gamma(n+\alpha +1)}\,.\]
For details and related material see \cite{Du3}, \cite{R1}, 
\cite{R-V}, \cite{Ros},
 and references cited there.
\end{enumerate}
\end{examples}

For later references, we next list some further  known properties of  $K$.

\begin{theorem} Let $z,w\in {\b C}^N$, $x,y\in\reel^N$, $\lambda\in{\b
    C}$, and $g\in W$. Then:
\begin{enumerate}
\item[\rm{(1)}]   $K(z,0)= 1\>$,
 $\> K(\lambda z,w) = K(z,\lambda w)\>$, $\>K(z,w)=K(w,z)\>$,
 $\>K(-ix,y)=\overline{K(ix,y)}\>$, and $K(g(x),g(y))=K(x,y)$; 
\item[\rm{(2)}]  $ T_j^x K(x,y) =  y_jK(x,y)$ for $j=1,\ldots,N$,
 where the superscript  $x$ indicates that the operators act with
 respect to the 
$x$-variable;
\item[\rm{(3)}]  For each $x\in\reel^N$ there is a unique probability measure $\mu_x\in M^1(\reel^N)$
 with \hfill\break 
$supp\>\mu_x\subset\{z\in\reel^N:\> |z|\le|x|\}$
   such that
$$K(x,y)=\int_{\reel^N} e^{\langle z,y\rangle}\> d\mu_x(z) \quad\quad{ for\>\>\> all}\>\>\>\> y\in{\b C}^N.$$
 In particular, $K(x,y)>0$ for all $x,y\in\reel^N$.
\item[\rm{(4)}] For all multi-indices $\nu\in\ganz_+^N$,
$$ |\partial_z^\nu K(x,z)|\le |x|^{|\nu|} \cdot e^{|x|\cdot |{\rm Re}\> z|}.$$
In particular,  $|K(z,w)| \le e^{|z||w|}$ and  $ |K(ix,y)|\le 1$.
\end{enumerate}
\end{theorem}  

\begin{proof} Parts  (1) and (2) follow from the properties of $K$ in Section 2.2;
 cf. [D2], [D3]. Part 
 (3) is shown in [R3], and Part (4) is a consequence of Part (3).
\end{proof}

The generalized exponential function $K$ gives rise to an integral transform, 
called the Dunkl transform on $\b R^N$, which was introduced in [D3]. 
 To emphasize the similarity to the classical Fourier transform, we use the following notion:

\begin{Dutr} The Dunkl transform  associated with $W$ and $k\geq 0$ is given by
$$\hat.: L^1(\b R^N, w_k(x)dx)\to C_b(\b R^N);\>\>\>
 \hat f(y):= \int_{\b R^N} f(x)\,K(-iy,x)\,w_k(x)dx \quad 
(y\in \b R^N).$$
 The Dunkl transform of a function $f\in L^1(\b R^N, w_k(x)dx)$
 satisfies $\|\hat f\,\|_\infty \le  \|f\|_{1,w_k(x)dx}$  by Theorem 2.4(4).
Moreover, according to  [D3], [dJ], and [R-V], many results from classical Fourier
analysis on $\reel^N$ have analogues for the Dunkl transform, like the
$L^1$-inversion theorem, the lemma of Riemann-Lebesgue, and 
Plancherel's formula.

We next extend the Dunkl transform to measures. We
 denote the Banach space of all $\b C$-valued bounded Borel measures on $\reel^N$ by
$M_b(\reel^N)$, while  $M^1(\reel^N)$ is the subspace 
consisting of all  probability measures on $\b R^N$. 

The Dunkl transform of a measure $\mu\in M_b(\reel^N)$
 is given by $\hat\mu(y) :=\int_{\reel^N} K(-iy,x)\> d\mu(x)$
$(y\in\reel^N$); it satisfies $\hat\mu\in C_b(\reel^N)$ with 
 $\|\hat \mu\|_\infty \le \|\mu\|$; cf. Theorem 2.4(4). Moreover:
\begin{enumerate}
\item[\rm{(1)}] If $\mu\in M_b(\reel^N)$ and $f\in L^1(\reel^N, w_k(x)dx)$, then
 $$\int_{\reel^N}\hat \mu(x)\> f(x)\> w_k(x) \> dx = \int_{\reel^N}\hat f \> d\mu;$$
\item[\rm{(2)}]  If $\mu\in M_b(\reel^N)$ satisfies $\hat\mu\equiv0$, then $\mu=0$.
\end{enumerate}
In fact,  Part (1) follows from Fubini's theorem, and Part (2) follows from Part (1) and the
fact that  $( L^1(\b R^N, w_k(x)dx))^\wedge$ 
is  $\|.\|_\infty$-dense in $C_0(\reel^N)$; see [dJ].
\end{Dutr}

We next turn to Dunkl's Laplace operator and the associated heat semigroup:

\begin{Lap} The generalized Laplacian $\Delta_k$ associated with 
 $W$ and  $k\geq 0$ is defined by
\begin{equation}
\Delta_k f \,:=\,\sum_{j=1}^N T_j^2 f\,;
\quad\quad f\in C^2(\reel^N).
\end{equation}
It is shown in [R2] that   $\Delta_k$ is a closable linear operator 
 on $C_0(\reel^N)$ and that its  closure
(again denoted by $\Delta_k$) generates a positive, strongly
 continuous contraction semigroup $(e^{t\Delta_k})_{t\ge0}$ on   
$C_0(\reel^N)$. This
 semigroup  is given  explicitly in terms of the following 
generalized heat kernels.
\end{Lap}

\begin{heat}
The generalized heat kernel $\Gamma_k$ is defined by
\begin{equation}
 \Gamma_k(x,y,t):=\, \frac{c_k}{(4t)^{\gamma +N/2}}\,e^{-(|x|^2 + |y|^2)/4t}\,
K\Bigl(\frac{x}{\sqrt{2t}},\frac{y}{\sqrt{2t}}\Bigr)
\quad\quad (x,y\in \b R^N,\> t>0)
\end{equation}
where  $c_k$ is given in (2.1). The  kernel
 $\Gamma_k$ has the following properties (see Lemma 4.5 in [R2]): Let 
$x,y,z\in\reel^N$ and $t>0.$ Then 
\begin{enumerate}
\item[\rm{(1)}] $\,\displaystyle  \Gamma_k(x,y,t) \,=\,  \Gamma_k(y,x,t) \,=\,
\frac{c_k^2}{4^{\gamma+N/2}} \int_{\b R^N} e^{-t|\xi|^2}\, K(ix,\xi)\, 
K(-iy,\xi) \, w_k(\xi) d\xi\,,$ and by the inversion formula for the
Dunkl transform,
$\Gamma_k(x,.,t)^\wedge(z)=e^{-t|z|^2}\cdot K(-ix,z)$.
\item[\rm{(2)}]  $\,\displaystyle \int_{\b R^N} \Gamma_k(x,y,t)\,w_k(x) 
dx\,=\,1\quad$ and  $\quad\displaystyle |\Gamma_k(x,y,t)|\,\leq\,
\frac{M_k}{t^{\gamma+N/2}}\,e^{-(|x|-|y|)^2/4t}.$  
\end{enumerate}
Moreover, the integral operators
\begin{equation}
 H(t)f(x):=\int_{\b R^N} \Gamma_k(x,y,t) f(y)\,w_k(y)dy\quad\quad
 {\rm for}\>\> t>0 \quad\quad{\rm and}\quad\quad H(0)f:=f
\end{equation}
 are related to the  semigroup
$(e^{t\Delta_k})_{t\ge0}$ by the fact that
 for all  $f\in C_0(\reel^N)\cup {\cal P} $
  \begin{equation}
e^{t\Delta_k}f=H(t)f \quad\>(t\geq 0).
\end{equation}
For  $f\in C_0(\reel^N)$ this is shown  in [R2].
 Moreover, by Proposition 2.1 of [D3] we have
\begin{equation}
p(x) =\int_{\reel^N} \Gamma_k(x,y,1/2)\> (e^{-\Delta_k/2} p)(y)
 \> w_k(y)\> dy\quad\quad{\rm for}\quad p\in{\cal P}.
\end{equation}
This proves (2.9) for $t=1/2$, as
   $e^{-\Delta_k/2}$ is the inverse  of  $e^{\Delta_k/2}$ on
 ${\cal P}$. The  general case $t>0$ follows by renormalization (see
   Lemma 2.1 of [R2]).
\end{heat}

We next turn to a probabilistic interpretation of the generalized heat kernels:

\begin{Gau} The $k$-Gaussian distribution $P_t^\Gamma(x,.)\in
  M^1(\reel^N)$ with "center" $x\in\reel^N$ and "variance parameter"
  $t>0$ is given by
\begin{equation}
P_t^\Gamma(x,A):= \int_A \Gamma_k(x,y,t)\> w_k(y)\> dy
\quad\quad(A\subset\reel^N\>\>\text{a Borel set}).
\end{equation}
In particular,
\[ P_t^\Gamma(0,A)\,=\, \frac{c_k}{(4t)^{\gamma+N/2}}\int_A
e^{-|y|^2/4t}\,w_k(y)dy.\]
It follows readily from the statements of Section 2.7 that the 
 $k$-Gaussian distributions $P_t^\Gamma(x,.)$ ($t>0$) have the
 following properties:
\begin{enumerate}
\item[\rm{(1)}] The Dunkl transforms of the probability measures 
$P_t^\Gamma(x,.)$ ($t\ge0,x\in\reel^N$)
satisfy
$$P_t^\Gamma(0,.)^\wedge(y) = e^{-t|y|^2}\quad{\rm and}\quad 
P_t^\Gamma(x,.)^\wedge(y)= K(-ix,y)\cdot P_t^\Gamma(0,.)^\wedge(y)
 \quad{\rm for}\quad  y\in\reel^N.$$
\item[\rm{(2)}] For each $t>0$, $P_t^\Gamma$ is a Markov kernel on
  $\reel^N$, and $(P_t^\Gamma)_{t\ge0}$ forms a semigroup of  Markov kernels on
  $\reel^N$; see Section 3.5 of  [R-V] for details.
\end{enumerate}
\end{Gau}

\section{ Moment functions}

In this section we first introduce homogeneous polynomials $m_\nu$
which generalize the monomials $x^\nu$ on $\reel^N$. These monomials
 will be called generalized moment functions and 
will be used to define generalized moments of $k$-Gaussian
distributions. 
Our approach is motivated by similar notions in [Bl-He] and references
there in the setting of probability measures on hypergroups. In the
sequel, a reflection group $W$ with root system $R$ and
multiplicity function $k\ge0$ is fixed.

\begin{Mom} As the Dunkl kernel $K$ is analytic on $\comp^{\,N\times N}$,
there exist unique analytic coefficient functions $m_\nu$ ($\nu\in \ganz_+^N$) on $\comp^N$ with
\begin{equation}
K(x,y)=\sum_{\nu\in\ganz_+^N}\frac{m_\nu(x)}{\nu!}y^\nu
\quad\quad (x,y\in\comp^N).
\end{equation}
The restriction of $m_\nu$ to $\reel^N$ is called the $\nu$-th moment
  function. It is  given explicitly  by
\begin{equation}
m_\nu(x)= (\partial_y^\nu K(x,y))|_{y=0}= V(x^\nu),
\end{equation}
where the first equation is clear by (3.1) and  the second one
follows from 
\[\partial_y^\nu K(x,y))|_{y=0}\,= \,\partial_y^\nu(V_x
e^{\langle x,y\rangle})|_{y=0} \,=\, V_x(\partial_y^\nu e^{\langle
  x,y\rangle}|_{y=0})=V(x^\nu);\]
see Section 2.2.
In particular, the homogeneity of $V$  ensures that
 $m_\nu \in \mathcal P_{|\nu|}$. 
Moreover, for each $n\in \b Z_+$, the moment functions $(m_\nu)_{|\nu|=n}$ 
 form a basis of the space ${\cal P}_n$.
\end{Mom} 

We have the following Taylor -type formula involving the moment
functions $m_\nu$:

\begin{proposition} Let $f:\comp^N\to\comp$ be analytic in
  a neighborhood of 0. Then
$$f(y)=\sum_{n=0}^\infty \sum_{ |\nu|= n}  \frac{m_\nu(y)}{\nu!}  \> 
T^\nu f(0),$$
where the series $\sum_{n=0}^\infty$ converges absolutely and uniformly in a neighborhood of 0.
\end{proposition}

\begin{proof} Assume first that $f\in{\cal P}$. As $V{\cal P}_n={\cal P}_n$, we have
$\partial^\nu f(0) = V\partial^\nu f(0) =T^\nu Vf(0)$. Thus,
$$f(y)= \sum_\nu \frac{y^\nu}{\nu!} \> T^\nu Vf(0) \quad\quad{\rm and}\quad\quad
(V^{-1}f)(y) =\sum_\nu \frac{y^\nu}{\nu!} \> T^\nu f(0),$$
which gives
$$f(y)=\sum_\nu \frac{m_\nu(y)}{\nu!}  \> T^\nu f(0).$$
The assertion now follows from the corresponding results for the classical case.
\end{proof}

Using the modified moment functions $m_\nu$, we next introduce
the modified moments of the $k$-Gaussian measures $P_t^\Gamma(x,.)$.

\begin{Mod} For $t>0,\> x\in\reel^N$ and $\nu\in\ganz_+^N,$   
the $\nu$-th  modified
  moment  of  $P_t^\Gamma(x,.)$ is defined by
\begin{equation}
m_\nu(P_t^\Gamma(x,.)):=\int_{\reel^N}m_\nu\> dP_t^\Gamma(x,.).
\end{equation}
These modified moments are closely related  to the classical
moments of the classical normal distributions on $\reel^N$:
\end{Mod}

\begin{lemma} Let $\nu\in\ganz_+^N$, $t>0$, and $x\in\reel^N$. Then:
\begin{enumerate}
\item[{\rm(1)}] $\displaystyle m_\nu(P_t^\Gamma(x,.))
= i^{|\nu|}\cdot \partial^\nu_y\,
  P_t^\Gamma(x,.)^\wedge(y)\Bigl|_{y=0}=
  i^{|\nu|}\cdot
\partial^\nu_y\Bigl( K(x,-iy)\cdot
e^{-t|y|^2}\Bigr)
 \Bigl|_{y=0}$;
\item[{\rm(2)}] $\displaystyle m_\nu(P_t^\Gamma(x,.))= \sum_{\rho\le\nu}
 \binom{\nu}{\rho} m_\rho(P_t^\Gamma(0,.))\cdot m_{\nu-\rho}(x)$;
\item[{\rm(3)}] $\displaystyle  m_\nu(P_t^\Gamma(0,.))
=  \begin{cases}
   \frac{(2\mu)!}{\mu!}\cdot t^{|\mu|}   &\text{if $\nu=2\mu$
     for some $\mu\in\ganz_+^N$}\\
   0  &\text{ otherwise.}
\end{cases}$
\end{enumerate}
 \end{lemma} 

\begin{proof} 
\begin{enumerate}
\item[(1)] The first equation is obtained from (3.2) and inductive use
 of the dominated convergence theorem (which is applicable by Theorem 2.4(4)).
The second assertion is clear.
\item[(2)] By Part (1), Eq.~(3.2) and the Leibniz rule for partial
 derivatives of products we obtain
\begin{align}
 m_\nu(P_t^\Gamma(x,.))\,=&\,
 i^{|\nu|}\cdot \sum_{\rho\in \b Z_+^N,
 \>\rho\leq \nu} \binom{\nu}{\rho}\,
  \partial_y^\rho\bigl(e^{-t|y|^2}\bigr)\Bigr|_{y=0}\cdot 
\partial_y^{\nu-\rho}K(x,-iy)\Bigr|_{y=0} \notag\\
 =\,& \sum_{\rho\leq\nu} \binom{\nu}{\rho}\,
  m_\rho(P_t^\Gamma(0,.))\cdot m_{\nu-\rho}(x).
\end{align}
\item[(3)] This follows from Part (1) for $x=0$ and the power series
  of the exponential function.
\end{enumerate}
\end{proof}

\section{Appell systems for $k$-Gaussian semigroups}

Based on the moment functions of the previous section and certain
generating functions,
 we now  construct two systems $(R_\nu)_{\nu\in\ganz_+^N}$ and  
$(S_\nu)_{\nu\in\ganz_+^N}$
of functions on $\reel\times\reel^N$ associated
 with the $k$-Gaussian semigroup $(P_t^\Gamma)_{t\ge0}$.
These systems, called Appell characters and cocharacters,
satisfy several algebraic relations, the most important being
the biorthogonality established in Theorem 4.5 below.
Among other results, we  present a new proof for  a generalization
of a formula of Macdonald [M] due to Dunkl [D2]. Our approach  and our
 notations are motivated by related concepts 
in non-Gaussian white-noise-analysis (see [ADKS], [Be-K]) and the
 umbral calculus in
 [Ro].

\begin{app} As the Dunkl kernel $K$ is analytic, we have a power
 series expansion of the form
\begin{equation}
\frac{K(x,-iy)}{P_t(0,.)^\wedge(y)}=
K(x,-iy)\cdot e^{t |y|^2}= \sum_{\nu\in\ganz_+^N}
 \frac{(-iy)^\nu}{\nu!} R_\nu(t,x) 
\quad\quad{\rm for}\quad t\ge0 \>\>\text{ and }\>x,\,y\in \b R^N,
\end{equation}
with certain functions $R_\nu$ on $[0,\infty[\times \b R^N$. Analogous
to  (3.4), these are
given by
\begin{equation}
 R_\nu(t,x)  =\,i^{|\nu|} \partial_y^\nu\Bigl(K(x,-iy)\cdot
 e^{t|y|^2}\Bigr)\Bigr|_{y=0}\, 
= \,\sum_{\rho\leq\nu} \binom{\nu}{\rho}  a_{\nu-\rho}(t)\cdot m_\rho (x),
 \end{equation}
 where
\begin{equation}
a_\lambda(t):= i^{|\lambda|}\cdot \partial_y^\lambda(
e^{t|y|^2})\bigr|_{y=0}
=  \begin{cases}
   \frac{(2\mu)!}{\mu!}\cdot (-t)^{|\mu|}   &\text{if $\lambda=2\mu$
     for $\mu\in\ganz_+^N$}\\
   0  &\text{ otherwise.}
\end{cases}
\end{equation}
In particular, the $R_\nu$ are real polynomials in the $(N+1)$
variables $(t,x)$ of degree $|\nu|$, and after analytic continuation,
$R_\nu(t,.)$ is a polynomial of degree $|\nu|$ for each $t\in\reel$. 
The polynomials  $R_\nu$  are called the
 {\it Appell characters}
associated with the semigroup $(P_t^\Gamma)_{t\ge0}$. 
\end{app}

We next collect some properties and examples of Appell characters.

\begin{lemma} In the setting of Section 4.1, the following  holds for all $\nu\in  \ganz_+^N$:
\begin{enumerate}
\item[{\rm(1)}] Inversion formula: 
 For all  $x\in\reel^N$ and $t\in\reel$,
$$m_\nu(x)=  \sum_{\rho\in\ganz_+^N,\> \rho\le \nu}
 \binom{\nu}{\rho} a_{\nu-\rho}(-t)\cdot R_\rho(t,x).$$
\item[{\rm(2)}] For all  $t\in\reel$ and $n\in \b Z_+$, the family
 $(R_\nu(t,.))_{\nu\in\ganz_+^N,\> |\nu|\le n}$
is a basis of the space $\bigoplus_{j=1}^n {\cal P}_j$ of all  polynomials of degree  at most $n$.
\item[{\rm(3)}] For  $x\in\reel^N$ and  $t\geq 0$, 
$$\int_{\reel^N} R_\nu(t,y)\> dP_t^\Gamma(x,.)(y)= m_\nu(x).$$
\item[{\rm(4)}] For $t>0$ and $x\in\reel$,  $R_\nu(t,x)= 
\sqrt{t}^{|\nu|}\cdot R_\nu(1,x/\sqrt t)$.
\item[{\rm(5)}] For all  $x\in\reel^N, \>t\in\reel$ and
  $j\in\{1,\ldots,N\}$, 
\[T_j R_{\nu+e_j}(t,x)= (\nu_j+1)\cdot R_\nu(t,x);\]
 here the Dunkl operator $T_j$ acts with
respect to the variable $x$ and $e_j\in\ganz_+^N$ is the $j$-th unit vector.
\end{enumerate}
\end{lemma}

\begin{proof}
\begin{enumerate}
\item[(1)]  Using (4.1) and (4.3), we obtain 
\begin{align} 
K(x,-iy) &=  e^{-t|y|^2}\cdot \Bigl(  e^{t |y|^2}K(x,-iy) \Bigr)
= \Bigl( \sum_{\lambda\in \ganz_+} \frac{ a_\lambda(-t)}{\lambda !}\>
 (-iy)^\lambda\Bigr)
\Bigl( \sum_{\rho\in \ganz_+} \frac{(-iy)^\rho}{\rho!} R_\rho(t,x)  \Bigr)\notag\\
&=\sum_{\nu\in \ganz_+}\Bigl( \sum_{\rho\le \nu} \binom{\nu}{\rho}
 a_{\nu-\rho}(-t) R_\rho(t,x)\Bigr)
\frac{(-iy)^\nu}{\nu!} \notag
\end{align}
A comparision of this expansion with Eq.~(3.1) leads to Part (1).
\item[(2)] This follows from Part (1) of this lemma,
 and the fact that $(m_\nu)_{|\nu|=j}$ is a basis of 
 ${\cal P}_j$.
\item[(3)] Comparison of formula (4.3)  and  Lemma 3.4(3) shows that  
$m_\lambda(P_t^\Gamma(0,.))
=a_\lambda(-t)$. Hence, by (4.2) and Lemma 3.4(2), 
\begin{align} 
\int_{\reel^N} R_\nu(t,y)\> dP_t^\Gamma(x,.)(y)\,&=
\,\sum_{\rho\le\nu}  \binom{\nu}{\rho}
 a_{\nu-\rho}(t) \cdot\int_{\reel^N} m_\rho(y) \>  dP_t^\Gamma(x,.)(y)\notag\\
&= \,\sum_{\rho\le\nu}  \binom{\nu}{\rho}
 a_{\nu-\rho}(t) \cdot\Bigl( \sum_{\lambda\le\rho}  \binom{\rho}{\lambda}
m_\lambda(P_t^\Gamma(0,.))\cdot m_{\rho-\lambda}(x)\Bigr)\notag\\
&=\,\sum_{\rho\le\nu}  \binom{\nu}{\rho}
 a_{\nu-\rho}(t) \cdot R_\rho(-t,x).\notag
\end{align}
The assertion now follows from the inversion formula of Part (1). 
\item[(4)] This is a consequence of  the homogeneity of the moment functions
  $m_\nu$ and of (4.2) and (4.3).
\item[(5)] Note first that by (3.2) and the intertwining property of $V$,
\[ T_j m_{\nu+e_j}\,=\, (\nu_j +1)\cdot m_\nu \quad  (j=1,\ldots,N\,,
\> \nu\in \b Z_+^N).\]
This, together with  identity (4.2)  and Proposition 3.2, yields
\begin{align}
T_j  R_{\nu+e_j}(t,x) &=\sum_{\rho\le \nu+e_j} \binom{\nu+e_j}{\rho} 
T_jm_\rho(x)\cdot  a_{\nu+e_j-\rho}(t) =
 \sum_{\rho\le \nu} \binom{\nu+e_j}{\rho+e_j}(\rho_j+1) m_\rho(x) \cdot  a_{\nu-\rho}(t)\notag\\
& =  (\nu_j+1)\cdot  \sum_{\rho\le \nu} \binom{\nu}{\rho}
m_\rho(x)\cdot a_{\nu-\rho}(t)
  =(\nu_j+1)\cdot R_\nu(t,x).
\notag
\end{align}
\end{enumerate}
\end{proof}

\begin{examples}\begin{enumerate}
\item[(1)]  In the classical case $k=0$ with  $m_\nu(x):= x^\nu$, Eq.~(4.3) leads to
\begin{equation}
R_\nu(t,x) = \sqrt t^{|\nu|}\cdot\widetilde H_\nu\Bigl(\frac{x}{2\sqrt t}\Bigr)
 \quad\quad\quad (x\in\reel^N,\> \nu\in\ganz_+^N,\>\>
 t\in\reel\setminus \{0\})
\end{equation}
where the $\widetilde H_\nu$ are the classical, $N$-variable Hermite polynomials defined by
$$\widetilde H_\nu(x) =\prod_{i=1}^N    H_{\nu_i}(x_i) \quad\quad{\rm with}\quad\quad
  H_{n}(y) = \sum_{j=0}^{\lfloor n/2\rfloor} \frac{(-1)^j\> n!}{j!\> (n-2j)!}\> (2y)^{n-2j};$$
c.f. Section 5.5 of [Sz] for the one-dimensional case.
\item[(2)] If $N=1$, $W=\b Z_2$  and $k\ge0$, then   Section
  2.3(2) easily leads to an explicit formula for 
the moment functions $m_n$. Using then (4.2),
  we finally obtain 
\begin{align}  R_{2n}(t,x) &=\, (-1)^n 2^{2n} n! \> t^n\>
 L_n^{(k-1/2)}(x^2/4t) \qquad\,\text{ and}\notag\\
R_{2n+1}(t,x) &=\,(-1)^n 2^{2n+1} n! \> t^n\>x\>
L_n^{(k+1/2)}(x^2/4t)\notag
\end{align}
for $n\in\ganz_+$; here the $L_n^{(\alpha)}$ are the  Laguerre
polynomials (see Section 5.1 of [Sz]) given by 
\[ L_n^{(\alpha)}(x) =\, \frac{1}{n!} x^{-\alpha}e^x\cdot\frac{d^n}{dx^n}
(x^{n+\alpha}e^{-x})\,=\, \sum_{j=0}^n \binom{n+\alpha}{n-j} 
\frac{(-x)^j}{j!}\,.\]
 The polynomials $(R_n)_{n\ge0}$ are called generalized Hermite
 polynomials and studied e.g. in [Ros].
For each $t>0$ the polynomials $(  R_{n}(t,.))_{n\ge0}$ are orthogonal
 with respect to the $k$-Gaussian measure
$$dP_t^\Gamma(0,.)(x) = \frac{\Gamma(k+1/2)}{(4t)^{k+1/2}}\>|x|^{2k}\>
e^{-x^2/4t}\>dx.$$
\end{enumerate}
\end{examples}

 An uninformed reader might suggest from these examples  that 
$k$-Gaussian Appell characters are always orthogonal with repect to
$P_t^\Gamma(0,.)$ for $t>0$. This is however not correct in general
for the  $S_N$- and $B_N$-cases; see Section 8 of [R-V].
To overcome this problem, we  introduce so-called Appell
cocharacters, which turn out to form  biorthogonal systems for the
Appell characters.

\begin{appco} 
 The    non-centered $k$-Gaussian measures  $P_t^\Gamma(x,.)$ 
admit 
$P_t^\Gamma(0,.)$-densities  $\theta_t(x,.)$ for $t>0,x\in\reel^N$. These
densities are given by
\begin{equation}
\theta_t(x,y):= \frac{dP_t^\Gamma(x,.)(y)}{dP_t^\Gamma(0,.)(y)} = 
e^{-|x|^2/4t} \> K(x,y/2t) =
\sum_{n=0}^\infty \sum_{|\nu|=n} \frac{m_\nu(x)}{\nu!} \> S_\nu(t,y)
\end{equation}
where, in view of Proposition 3.2, the coefficients $S_\nu$ are given by
$$S_\nu(t,y) = T_x^\nu\bigl( e^{-|x|^2/4t} \> K(x,y/2t)\bigr)\bigl|_{x=0}.$$
We mention that the function $\theta_t$ (with $t=1$) appears also as
the generating function of the generalized Hermite polynomials
associated with $W$ and $k$ in [R2]. By Proposition 3.8 of [R2], 
 the convergence of the series $\sum_{n=0}^\infty$ in (4.5)
 is locally uniform on $\comp^N\times\comp^N$.
The functions $S_\nu$ are called the {\it Appell cocharacters} of the
 $k$-Gaussian semigroup $(P_t^\Gamma)_{t\ge0}$. 

Using the homogeneity of $m_\nu$, we obtain the following analogue of
Lemma 4.2(4):
\begin{equation}
S_\nu(t,y)= \Bigl(\frac{1}{\sqrt t}\Bigr)^{|\nu|}\cdot S_\nu(1,
y/\sqrt t\>) \quad\quad (t>0).
\end{equation}
\end{appco}

As announced, Appell   characters and cocharacters have  the following
 biorthogonality property:

\begin{theorem}
 Let $t>0$, $\nu,\rho\in\ganz_+^N$, and $p\in {\cal P}$ a polynomial
 of degree less than $|\nu|$. Then:
\begin{enumerate}
\item[{\rm(1)}]  $\displaystyle \int_{\reel^N} R_\nu(t,y) \,
  S_\rho(t,y)\> dP_t^\Gamma(0,.)(y)
 = \nu!\> \delta_{\nu,\,\rho}\,;$
\item[{\rm(2)}]  $\displaystyle \int_{\reel^N} p(y)\, S_\nu(t,y)\> 
dP_t^\Gamma(0,.)(y)=
\int_{\reel^N} p(y)\, R_\nu(t,y)\> dP_t^\Gamma(0,.)(y)=0$.
\end{enumerate}\end{theorem}

\begin{proof}  We use the definition of $\theta_t$ and Lemma 4.2(3) 
and conclude that for $x\in\reel^N$,
\begin{align}
m_\nu(x) &= \int_{\reel^N} R_\nu(t,y) \,\theta_t(x,y)\> \>
dP_t^\Gamma(0,.)(y)
= \int_{\reel^N} \sum_{n=0}^\infty\sum_{|\rho|=n} R_\nu(t,y)
 S_\rho(t,y)\> \frac{m_\rho(x)}{\rho!}  \>  dP_t^\Gamma(0,.)(y)\notag\\
 &= \sum_{n=0}^\infty\sum_{|\rho|=n} \frac{m_\rho(x)}{\rho!}  \>
 \int_{\reel^N} R_\nu(t,y)\>  S_\rho(t,y)\> dP_t^\Gamma(0,.)(y),
\end{align}
where we still have to justify that the summation can be interchanged
with 
the integration. 
For this, we restrict our attention to the case $t=1/4$, as the
general case then follows by renormalization (see Lemma 4.2(4) and
formula (4.6).)  We follow the proof of Proposition
3.8 of [R2] and decompose $\theta_{1/4}(x,y)$ into its
 $x$-homogeneous parts:
$$\theta_{1/4}(x,y)= \sum_{n=0}^\infty L_n(y,x) \quad\quad{\rm with} \quad\quad
L_n(y,x)= \sum_{|\rho|=n} \frac{m_\rho(x)}{\rho!} S_\rho(1/4\,,y).$$
The estimations of Theorem 2.4(4) imply that
$$|L_{2n}(y,x)|\le \frac{|x|^{2n}}{n!}\cdot (1+2|y|^2)^{n} \qquad(n\in
\b Z_+),$$
and a similar estimations holds if the degree is odd; for details see the proof of 3.8 in [R2].
Therefore,
$$\sum_{n=0}^\infty\int_{\reel^N} 
 |L_n(y,x)|\> R_\nu(1/4,y)\> dP_{1/4}^\Gamma(0,.)(y) <\infty.$$
The dominated convergence theorem now justifies the last step in
(4.7) for $t=1/4$.
 Part (1) of the theorem  is now a immediate consequence of relation
 (4.7). Finally, Part
 (2) follows from  Part (1), Lemma 4.2(2), and Section 4.4.
\end{proof}

Together with Lemma 4.2(2), Theorem 4.5(1) in particular implies that 
for each $n\in \b Z_+$, the family $(R_\nu(t,.))_{\nu\in \b Z_+,\,
  |\nu|\leq n}$ is a basis of $\bigoplus_{j=1}^n {\cal P}_j$.
The following result also  reflects the dual nature of Appell
characters and 
cocharacters.

\begin{proposition} Let $t\in\reel,\> x\in\reel^N$ and   $\nu\in\ganz_+^N$. Then 
$$R_\nu(t,x)= e^{-t\Delta_k}m_\nu(x)
 \quad\quad{ and} \quad\quad
S_\nu(t,x)= \bigl(\frac{1}{2t}\bigr)^{|\nu|} \cdot  e^{-t\Delta_k}x^\nu.$$
\end{proposition} 

\begin{proof}
In view of Section 2.7, Lemma 4.2(3) just says that  
$e^{t\Delta_k}R_\nu(t,x)=m_\nu(x)$ for $t\ge0$.
This shows the first part  for  $t\ge0$. As both sides are polynomials
in $t$ (write $e^{-t\Delta_k}$ as terminating power series of
operators acting on polynomials!), the first equation  holds generally.
 Let $\Delta_k^y$ be the $k$-Laplacian acting on the variable $y$, 
and let $V_x$ be the intertwining operator
acting on the variable $x$. Then
$$ e^{t\Delta_k^y}\Bigl( e^{-|x|^2/4t}\> K\Bigl(x,
\frac{y}{2t}\Bigr)\Bigr) = 
e^{-|x|^2/4t}\>\cdot e^{|x|^2/4t}\> K\Bigl(x, \frac{y}{2t}\Bigr)
=  K\Bigl(x, \frac{y}{2t}\Bigr) = V_x(e^{\langle x,y/2t\rangle}).$$
Now consider on both sides the homogeneous part $W_n$  of degree $n$ in the variable $x$.
 Using the left hand side, we obtain from (4.5) that
$$W_n= e^{t\Delta_k^y}\Bigl(\sum_{|\nu|=n} \frac{m_\nu(x)}{\nu!}\>  S_\nu(y,t)\Bigr)=
 \sum_{|\nu|=n} \frac{m_\nu(x)}{\nu!}\>  e^{t\Delta_k^y} \> S_\nu(y,t).$$
Using the right hand side, we conclude from the identity $V_x(x^\nu)=m_\nu(x)$ that
$$W_n= V_x\Bigl(\sum_{|\nu|=n} \frac{x^\nu}{\nu!}\>  (y/2t)^\nu\Bigr) = 
\sum_{|\nu|=n} \frac{m_\nu(x)}{\nu!}\>  (y/2t)^\nu.$$
A comparision of the corresponding coefficients leads to the second equation.
\end{proof}

We now combine Theorem 4.5 and Proposition 4.6 to rediscover 
 a   generalization of a formula of Macdonald [M] due to Dunkl [D2].
However,  our proof is  completely different from [D2].
We need the following notation: For a multiplicity function $k\ge0$, we introduce the bilinear form
\begin{equation}
[p,q]_k := (p(T)q)(0)
 \quad\quad\quad\quad{\rm for}\quad p,q\in{\cal P}.
\end{equation}

\begin{corollary} For all  $p,q\in{\cal P}$ and $t>0$,
$$[p,q]_k = \frac{1}{(2t)^{|\nu|}} \int_{\reel^N}
 e^{-t\Delta_k}(p) \, 
 e^{-t\Delta_k}(q)\> dP_t^\Gamma(0,.).$$
In particular, $[.,.]_k$ is a scalar product on ${\cal P}$.
\end{corollary}

\begin{proof}
Let $t>0$ and $\nu,\rho\in\ganz_+^N$. Then, by Theorem 4.5(1) and 
Proposition 4.6,
$$  \frac{1}{(2t)^{|\nu|}} \int_{\reel^N} e^{-t\Delta_k} 
(x^\nu)\, 
 e^{-t\Delta_k}(m_\rho)\> dP_t^\Gamma(0,.) = \nu!\cdot 
\delta_{\nu,\rho}.$$
On the other hand, as $V$ acts on ${\cal P}$ in a homogeneous way,
\begin{equation}
[x^\nu,m_\rho]_k = (T^\nu m_\rho)(0) = (T^\nu Vx^\rho)(0)= 
(V\partial^\nu x^\rho)(0) = 
 \nu!\cdot \delta_{\nu,\rho}\,.
\end{equation}
This yields the first statement. The second statement is clear.
\end{proof}

We give a further application of Theorem 4.5 for $t=1/2$. For this, we use the adjoint operator
$T_j^*$ of the Dunkl operator $T_j$ ($j=1,\ldots,N$) in $L^2(\reel^N,dP_{1/2}^\Gamma(0,.)),$ which is
given by
\begin{equation}
T_j^*f(x) =x_jf(x)-T_jf(x) = -e^{|x|^2/2}\cdot T_j\Bigl( e^{-|x|^2/2}\>f(x)\Bigr) \quad\quad(f\in{\cal P});
\end{equation}
see Lemma 3.7 of [D2]. (The second equation is a  consequence of the
product rule 2.2(3).)

\begin{corollary} For all  $\nu\in\ganz_+^N, \> j=1,\ldots,N, \>x\in\reel^N$ and $t>0,$
\begin{enumerate}
\item[(1)] $S_{\nu+e_j}(1/2, x) = T_j^*S_\nu(1/2, x)$;
\item[(2)] Rodriguez formula:
 $S_\nu(t, x)= (-1)^{|\nu|} e^{|x|^2/4t}\cdot T^\nu\bigl( e^{-|x|^2/4t}\bigr)$.
\end{enumerate}\end{corollary}

\begin{proof} 
  Using Theorem 4.5(1) and Lemma 4.2(5), we obtain for all
  $\rho\in\ganz_+^N$ that
\begin{align}
 \int_{\b R^N} R_{\rho+e_j}\cdot  T_j^*S_\nu\> dP^\Gamma_{1/2}(0,.)& =\,
  \int_{\b R^N} T_j R_{\rho+e_j}\cdot S_\nu\> dP^\Gamma_{1/2}(0,.)
=\,(\rho_j+1)\int_{\b R^N}  R_{\rho}\cdot S_\nu\>
  dP^\Gamma_{1/2}(0,.)\notag\\
& =\,
 \delta_{\rho,\nu}\cdot (\rho+e_j)! =\, 
\int_{\b R^N}  R_{\rho+e_j}\cdot S_{\nu+e_j}\> dP^\Gamma_{1/2}(0,.).\notag
\end{align}
As $\cal P$ is dense in  $L^2(\reel^N,dP_{1/2}^\Gamma(0,.))$, this
 implies Part (1).
 Part (2) for $t=1/2$ now follows from
(4.10), and the general case is a consequence of formula (4.6).
\end{proof}

\begin{remarks} 
 In the examples of Section 4.3 (i.e., in the classical case
  $k=0$ and in the one-dimensional case), Proposition 4.6 implies 
 that the systems $(S_\nu(t,.))_{\nu\in \b
  Z_+^N}$ and $(R_\nu(t,.))_{\nu\in \b Z_+^N}$ coincide -- up to
  multiplicative factors which depend on $t$ and $\nu$. Recall also
  that in these cases they  are orthogonal
  with respect to $P_{t}^\Gamma(0,.)$, and can be considered as
  generalized Hermite polynomials. 
In general, however, both systems fail to be orthogonal.
In any case, Corollary 4.7 allows to introduce orthogonal polynomials 
 with respect to $P_t^\Gamma(0,.)$. For this,
one has to choose an orthogonal basis
$(\phi_\nu)_{\nu\in\ganz_+^N}\subset{\cal P}$ with respect to the
scalar product $[.,.]_k$ in (4.8) with $\phi_\nu\in {\cal P}_{|\nu|}$
for $\nu\in\ganz_+^N$. (Note that by (4.9), $\mathcal P_n \,\bot\,
  \mathcal P_m$ for $n\not=m$).
It is then clear
  from Corollary 4.7 that
$(H_\nu:= e^{-t\Delta_k}\phi_\nu)_{\nu\in\ganz_+^N}$ forms a system
of orthogonal polynomials with respect to $P_{t}^\Gamma(0,.)$. These
 generalized Hermite polynomials are studied (for $t=1/4$) in [R2], and their
 relations to the Appell systems $(R_\nu)_{\nu\in\ganz_+^N}$ and 
$(S_\nu)_{\nu\in\ganz_+^N}$  are studied in [R-V]. We also refer to
 related investigations in [B-F, vD] and references given there.
\end{remarks}

\end{document}